\documentclass[12pt]{article}
\usepackage{graphicx,epsfig}
\usepackage{color}
\usepackage{amsmath,amssymb,mathrsfs,amsthm}
\usepackage{amsfonts}
\usepackage{multicol}
 \usepackage{graphicx}
\usepackage{color}
\usepackage{pdfcolmk}
\usepackage[nice]{nicefrac}
\usepackage{amsthm}
\usepackage{latexsym}
\usepackage{dsfont}
\usepackage{setspace}
\usepackage[a4paper, tmargin=1.5in, bmargin=1in, lmargin=1in, rmargin=1in, headheight=13.6pt]{geometry}
\usepackage[sort&compress,numbers]{natbib}
\usepackage[font=small,labelfont=bf,labelsep=space]{caption}
\makeatletter

\date{}

\def\@citex[#1]#2{\if@filesw\immediate\write\@auxout{\string\citation{#2}}\fi
  \def\@citea{}\@cite{\@for\@citeb:=#2\do
    {\@citea\def\@citea{,\linebreak[0]\hskip0pt plus .2em}%
      \@ifundefined{b@\@citeb}%
    {{\bf ?}\@warning{Citation `\@citeb' on page \thepage\space undefined}}%
      \hbox{\csname b@\@citeb\endcsname}}}{#1}}

\newtheorem{example}{Example}[section]

\newtheorem{rule-def}[theorem]{Rule}
\numberwithin{equation}{section}
\makeatletter

\begin{document}
\title{Adomian decomposition method for solving derivative-dependent doubly singular boundary value problems}
\author{Randhir Singh \thanks{Corresponding author. E-mail:randhir.math@gmail.com}\\
$\rm $\small {Department of Mathematics}\\
\small {Birla Institute of Technology Mesra,
 Ranchi-835215, India}} \maketitle{}
\begin{abstract}
 \noindent In this work, we apply Adomian decomposition method for solving nonlinear derivative-dependent doubly singular boundary value problems: $(py')'= qf(x,y,y')$. This method is based on the  modification of  ADM and new two-fold integral operator. The approximate  solution is obtained in the form of series with easily determinable components. The effectiveness  of the proposed approach is examined  by considering three examples and numerical results are compared with known results.
\end{abstract}
\section{Introduction}
Singular boundary value problems for ordinary differential equations arise very
frequently in many branches of applied mathematics and physics such as gas dynamics, chemical reactions, nuclear physics,  atomic structures, atomic calculations, and study of positive radial solutions of nonlinear elliptic equations and physiological studies  \cite{lin1976oxygen,chandrasekhar1939an,keller1955electrohydrodynamics,fermi1927metodo}.  Consider the following class of nonlinear  singular boundary value problems
 \begin{eqnarray}\label{sec1:1}
 (p(x)y'(x))'= q(x)f(x,y(x),y'(x)),~~0<x\leq 1,
\end{eqnarray}
with boundary conditions
\begin{align}\label{sec1:0}
y(0)=\eta_1,~\alpha_1y(1)+\beta_1y'(1)=\gamma_1,
\end{align}
where $\alpha_1>0$, $\beta_1\geq0$, $\gamma_1$ and  $\eta_1$ are any finite constants. The condition $p(0)=0$ characterizes that the problem (\ref{sec1:1}) is singular and if $q(x)$ is
allowed to be discontinuous at $x=0$ then the problem (\ref{sec1:1}) is  called doubly singular
\cite{Bobisud1990}. 

The study of singular boundary value problems  has attracted the attention of several researchers \cite{Bobisud1990,dunninger1986existence}. In \cite{dunninger1986existence}, for continuous $q(x)$ and in \cite{Bobisud1990}, for $q(x)$ not necessarily continuous at $x=0$ the existence as well as uniqueness of solution of equation \eqref{sec1:1} was discussed. The main difficulty of problem \eqref{sec1:1} is that  the singularity behavior  occurs at $x=0$. The establish of existence and uniqueness of solution of the singular problems \eqref{sec1:1}-\eqref{sec1:0} is presented in \cite{pandey2011constructive}.

Finding the solution of nonlinear SBVP \eqref{sec1:1}, most of  the methods designed for nonsingular boundary value problems suffer from a loss of accuracy or may even fail to converge \cite{ascher1995numerical}, because of the singular behavior  at $x=0$.  However, the finite difference method can be used to solve linear singular two point BVPs, but it may be difficult to solve nonlinear singular two-point boundary value problems. Moreover, the finite difference method requires some major modifications that include the use of some root-finding technique while solving nonlinear singular two-point BVPs.

In \cite{mohanty2005family}, the author has discussed existence and uniqueness of solutions of singular equation  $y''=f(x,y,y')$ and presented variable mesh methods for numerically solving such problems.  The higher order finite difference and cubic spline methods are carried out in \cite{ravi2004higher,ravi2005cubic} for the singular BVP $y''+(\nicefrac{k}{x}) y'=r(x)-q(x)y$, where $k$ is any constant.  Recently, some newly developed approximate methods have also been employed to solve special case of  SBVPs \eqref{sec1:1}. For example, Adomian decomposition method \cite{m2003decomposition,khuri2010novel,ebaid2011new}, the Homotopy perturbation method
 \cite{ramos2008series}, homotopy analysis method \cite{bataineh2009homotopy}  and variational iteration method \cite{das2016algorithm,singh2017optimal}.

\subsection{Standard Adomian decomposition method }
Adomian decomposition method has  received a lot of attention because it allows solution of both linear and nonlinear  ordinary differential, partial differential and integral equations. Adomian \cite{Adomian1994a} asserted that the ADM provides an efficient and computationally worthy method for generating approximate series solution for the large class of functional equations.
According to standard ADM    equation \eqref{sec1:1} can be rewritten as
\begin{align}\label{sec1:2}
Ly=Ny,
\end{align}
where $L\equiv\frac{d^2}{dx^2}$ is  linear  derivative operator;
 $Ny=-\nicefrac{p'}{p}y'+\nicefrac{q}{p}f(x,y,y')$ represents the nonlinear term.
The inverse  operator of  $L$ is defined  as
\begin{eqnarray}
L^{-1}(\cdot)=\int\limits_{0}^{x}\int\limits_{0}^{x}(\cdot) dx dx.
\end{eqnarray}
Operating the inverse linear operator  $L^{-1}(\cdot)$ on  both the sides
of \eqref{sec1:2}, we obtain
\begin{align}\label{sec1:3}
 y=y(0)+xy'(0)+L^{-1}Ny
\end{align}
Next, we  decompose the solution $y$ and the nonlinear function $Ny$  by an
infinite series as
\begin{align}\label{sec1:4}
y=\sum_{n=0}^{\infty}y_n,~~Ny=\sum_{n=0}^{\infty}A_n
\end{align}
where $A_n$ are Adomian polynomials that can be constructed for various classes of nonlinear functions  with the  formula given by Adomian and Rach \cite{adomian1983inversion}
\begin{eqnarray}\label{sec1:5}
A_n=\displaystyle\frac{1}{n!}\frac{d^n}{d
\lambda^n}\left[N\left(\sum_{k=0}^{\infty} y_k\lambda^k
\right)\right]_{\lambda=0},~~~n=0,1,2,...
\end{eqnarray}
Substituting  the series \eqref{sec1:4} into \eqref{sec1:3}, we obtain
\begin{align}\label{sec1:6}
\sum_{n=0}^{\infty}y_n=y(0)+xy'(0)+L^{-1}\sum_{n=0}^{\infty}A_n.
\end{align}
From the equation \eqref{sec1:6}, the ADM admits the following recursive  scheme
\begin{eqnarray}\label{e0}
\left.
  \begin{array}{ll}
  y_0=y(0)+x y'(0),\\
 y_{k+1}= L^{-1}A_{k},~~~~k\geq0.
\end{array}
\right\}
\end{eqnarray}
that will lead to the complete determination of components $y_n$ and hence the $n$-term truncated  is abstained as
\begin{eqnarray}
\psi_n(x)=\sum_{m=0}^{n-1} y_m(x),
\end{eqnarray}
The ADM has been applied to solve nonlinear boundary value problems for ordinary differential equations by many researchers   \cite{adomian1993,adomian1994,wazwaz2000approximate,
  wazwaz2001reliable,m2003decomposition,benabidallah2004application,jang2008two,khuri2010novel,ebaid2011new,singh2013solving}.
Solving nonlinear boundary value problems using standard ADM or modified ADM is always a computationally involved task as it requires the computation of undetermined coefficients in a sequence of nonlinear algebraic equations  which increases the computational work, (see \cite{benabidallah2004application,m2003decomposition,wazwaz2000approximate,wazwaz2001reliable}).

In order to avoid solving such nonlinear algebraic equations for nonlinear two-point  boundary value problems with derivative dependent source function, we use ADM which does not require any addition computational work for unknown constant based on the work \cite{singh2013solving,singh2014efficient,singh2013numerical,singh2014approximate,singh2014adomian,singh2016efficient}.

\section{Modified Adomian decomposition  method}
In this section, we apply improved ADM based on two-fold integral operator  for solving   nonlinear singular boundary value  problems. We  rewrite equation  \eqref{sec1:1} as:
\begin{align}\label{sec2:1}
Ly(x)=q(x)f(x,y,y'),
\end{align}
where $Ly(x)\equiv (p(x)y'(x))'$ is the linear differential operator to be inverted. Two-fold integral operator  $ L^{-1}(\cdot)$ regarded as the inverse operator of $L(.)$ is  proposed as
\begin{eqnarray}
L^{-1}(\cdot)\equiv\int\limits_{0}^{x}\frac{1}{p(s)}\int\limits_{s}^{1}(\cdot)dxds.
\end{eqnarray}
We operate $L^{-1}(\cdot)$ on  the left hand side of  \eqref{sec2:1} and use  $y(0)=\eta_1$ yields
\begin{align}\label{sec2:2}
\nonumber L^{-1}[(p(x)y'(x))']&=\int\limits_{0}^{x}\frac{1}{p(s)}\int\limits_{s}^{1}(p(x)y'(x))'dx
ds=\int\limits_{0}^{x}\frac{1}{p(s)}(p(1)y'(1)-p(s)y'(s))ds,\\
&=c_1\int\limits_{0}^{x}\frac{1}{p(s)}ds-y(x)+\eta_1=c_1\int\limits_{0}^{x}\frac{1}{p(s)}ds-y(x)+\eta_1.
\end{align}
We  again operate  $L^{-1}(.)$ on both sides of  \eqref{sec2:1} and use \eqref{sec2:2},  gives
\begin{eqnarray}\label{sec2:3}
y(x)=\eta+c_1\int\limits_{0}^{x}\frac{1}{p(s)}ds-\int\limits_{0}^{x}\frac{1}{p(s)}\int\limits_{s}^{1}q(x)f(x,y,y')dxds.
\end{eqnarray}
For simplicity, we set $$h(x)= \int\limits_{0}^{x}\frac{1}{p(s)}ds,~~
[L^{-1}(.)]_{x=1}=\int\limits_{0}^{1}\frac{1}{p(s)}\int\limits_{s}^{1}(.)dxds,~~~~~h'(1)=\frac{1}{p(1)}.$$
Then the equation \eqref{sec2:3} may be written as
\begin{eqnarray}\label{sec2:4}
y(x)=\eta_1+c_{1}h -[L^{-1}q f(x,y,y')].
\end{eqnarray}
 To eliminate  $c_1$ from \eqref{sec2:4},  we impose $\alpha_1y(1)+\beta_1y'(1)=\gamma_1$, and we obtain
 \begin{eqnarray}\label{sec2:5}
c_1=\frac{(\gamma_1-\alpha_1\eta)}{\alpha_1h(1)+\beta_1 h'(1)}+\frac{\alpha_1}{\alpha_1h(1)+\beta_1 h'(1)}
\big[L^{-1}q f(x,y,y')\big]_{x=1}.
\end{eqnarray}
Substituting the value of $c_1$ into \eqref{sec2:4}, we obtain
\begin{align}\label{sec2:6}
\displaystyle y(x)=\eta_1+\frac{(\gamma_1-\alpha_1\eta_1) h}{\alpha_1h(1)+\beta_1 h'(1)}+\frac{\alpha_1 h}{\alpha_1h(1)+\beta_1 h'(1)}
\bigg[L^{-1}q f(x,y,y')\bigg]_{x=1}-\big[L^{-1}q f(x,y,y')\big].
\end{align}
Substituting the series defined in \eqref{sec1:4} into  \eqref{sec2:6} gives
\begin{equation}\label{sec2:eq9}
\sum_{n=0}^{\infty}y_n=\eta_1+\frac{(\gamma_1-\alpha_1\eta_1)h}{\alpha_1h(1)+\beta_1 h'(1)}+\frac{\alpha_1 h}{\alpha_1h(1)+\beta_1 h'(1)}
\big[L^{-1}q\sum_{n=0}^{\infty} A_n\big]_{x=1}-\bigg[L^{-1}q\sum_{n=0}^{\infty} A_n \bigg].
\end{equation}
Comparing  the both sides of  \eqref{sec2:eq9}, we have
\begin{eqnarray}\label{f2}
\left.
  \begin{array}{ll}
\displaystyle    y_0=\eta_1,\\
\displaystyle y_1=\frac{(\gamma_1-\alpha_1\eta)h}{\alpha_1h(1)+\beta_1 h'(1)}+\frac{\alpha_1h}{\alpha_1h(1)+\beta_1 h'(1)}\big[L^{-1}qA_{0}\big]_{x=1}-\big[L^{-1}qA_{0} \big],\\
\vdots\\
\displaystyle y_{n+1}=\frac{\alpha_1h(x)}{\alpha_1h(1)+b
h'(1)}\big[L^{-1}q A_{n}\big]_{x=1}-\big[L^{-1}q A_{n}\big],~~~~n\geq1.
  \end{array}
\right\}
\end{eqnarray}
The  recursive scheme  \eqref{f2}
gives the complete determination of solution components $y_n$ of solution $y$ and hence the approximate  series solution $\psi_n$ can be obtained as
\begin{eqnarray}\label{sec2:eq10}
\psi_n=\sum_{m=0}^{n-1} y_m.
\end{eqnarray}

\section{Numerical illustrations and discussions}
In order to demonstrate the effectiveness and efficaciously of the proposed method, we have considered three complicated nonlinear singular examples. All the numerical results obtained by proposed ADM are compared with known results. Now,  we denote error functions as $E_n(x)=|\psi_n(x)-y(x)|$  and the maximum absolute errors as
 \begin{align}
E^{n}=\max_{0<x\leq1}E_n(x).
\end{align}

\begin{example}\label{sec3:ex1}
{\rm Consider the nonlinear  boundary value problem
\begin{eqnarray}\label{prob_1}
\left.
  \begin{array}{ll}
(x^{\alpha}y')'= \beta x^{\alpha+\beta-2}e^{y}\big( -x y'-\alpha-\beta+1 ),~~~~0<x\leq 1,\\
y(0)=\ln\frac{1}{4}~~~y(1)=\ln\frac{1}{5},
\end{array}
\right\}
\end{eqnarray}
with exact solution $y(x)=\ln(\frac{1}{4+x^{\beta}})$.}
\end{example}
Applying \eqref{f2} to the equation \eqref{prob_1}, with  $\alpha_1=1$, $\beta_1=0$,
 $\gamma_1=-\ln(5)$ and $\eta=-\ln(4)$, we have
\begin{align}\label{sec3:s1}
\left.
  \begin{array}{ll}
 y_0=-\ln(4),\\
 y_1=(\ln4-\ln5)x^{1-\alpha}+x^{1-\alpha}\int\limits_{0}^{1}x^{-\alpha}\int\limits_{x}^{1}x^{\alpha +\beta -2}A_{0}dxdx -\int\limits_{0}^{x}x^{-\alpha}\int\limits_{x}^{1}x^{\alpha +\beta -2}A_{0}dxdx,\\
y_{n+1}=x^{1-\alpha}\int\limits_{0}^{1}x^{-\alpha}\int\limits_{x}^{1}x^{\alpha +\beta -2}A_{n}dxdx-
\int\limits_{0}^{x}x^{-\alpha}\int\limits_{x}^{1}x^{\alpha +\beta -2}A_{n}dxdx,~~n\geq1.
  \end{array}
\right\}
\end{align}
The Adomian polynomials  for  $f(x,y,y')=-\beta \big(xe^{y}y'+e^{y}(\alpha+\beta-1) )$ with $y_0=-\ln(4)$ are obtained as:
\begin{align}\label{sec3:a1}
\left.
 \begin{array}{ll}
 A_0=-\beta e^{y_0}\big(x y'_0+(\alpha+\beta-1)\big)  \\
 A_1=-\beta e^{y_0}\big(x y'_1+y_1(\alpha +\beta -1) \big),\\
 A_2=-\beta e^{y_0}\big(x y'_2+y_1 y'_1)+( y_2+\nicefrac{y_1^2}{2} )(\alpha -\beta +1)\big),\\
 A_3=-\beta  e^y\big(x \left( y'_3+y'_2y_1+y'_1\left(\nicefrac{y_1^2}{2}+y_2\right)\right)+\left( y_3+y_1 y_2 +\nicefrac{y_1^3}{6} \right)(\alpha +\beta -1) \big),\\
\ldots
\end{array}
\right\}
\end{align}

For~$\alpha=0.5$,~$\beta=1$: Using \eqref{sec3:s1} and \eqref{sec3:a1},  we obtain
\begin{align*}
y_0&=-1.38629436,\\
y_1&=0.0268564x^{0.5}-0.25 x,\\
y_2&=-0.0267739x^{0.5}-0.00447607 x^{1.5}+0.03125 x^{2},\\
y_3&=-0.0003279x^{0.5}+0.0044623x^{1.5}-0.000045079x^{2}+0.00111902x^{2.5}-0.0052083x^{3},\\
y_4&=0.00025327x^{0.5}+0.0000546545 x^{1.5}+0.0000898816 x^{2}-0.0011159 x^{2.5}\\&+0.0000212874 x^{3}-0.0002797x^{3.5}+0.000976563 x^{4},\\
\ldots
\end{align*}
For~$\alpha=0.5$,~$\beta=3.5$: Using \eqref{sec3:s1} and \eqref{sec3:a1},  we obtain
\begin{align*}
y_0&=-1.38629436,\\
y_1&=0.+0.0268564 x^{0.5}-0.25 x^{3.5},\\
y_2&=0.-0.0253752 x^{0.5}-0.00587485 x^{4}+0.03125 x^{7},\\
y_3&=-0.00174107 x^{0.5}+0.00555081x^{4}-0.000070123x^{4.5}+0.00146871x^{7.5}-0.00520833x^{10.5},\\
y_4&=0.000230726 x^{0.5}+0.000380859 x^{4}+0.000132511 x^{4.5}-5.6497931825\times 10^{-7} x^{5}\\&-0.0013877 x^{7.5}+0.0000347878 x^{8}-0.000367178 x^{11}+0.000976563 x^{14},\\
\ldots
\end{align*}
The maximum absolute error $E^{(n)}$ for $n=5,8$, and $10$ are listed in Table \ref{tab1} and \ref{tab2} for different values of $\alpha$ and $\beta$.
\begin{table}[htbp]\caption{Maximum error of Example~\ref{sec3:ex1}, when $\beta=1$.}
\centering
\vspace{-.3cm}
\renewcommand{\arraystretch}{1.0}
\setlength{\tabcolsep}{0.25in}
\begin{tabular}{c|c cc}
\cline{1-4}
$\alpha$ &  $E^{(5)}$               &$E^{(8)}$               & $E^{(10)}$ \\
\cline{1-4}
0.25    &1.11205$\times 10^{-5}$  &2.30301$\times 10^{-8}$   &3.54171$\times 10^{-10}$ \\
0.5     &1.52386$\times 10^{-5}$  &2.68725$\times 10^{-8}$   &6.11240$\times 10^{-10}$ \\
0.75    &1.62907$\times 10^{-5}$  &4.91114$\times 10^{-8}$   &9.22449$\times 10^{-10}$ \\
\hline
\end{tabular}
\label{tab1}
\end{table}
\begin{table}[htbp]\caption{Maximum error of Example \ref{sec3:ex1}, when $\beta=3.5$.}
\centering
\vspace{-.3cm}
\renewcommand{\arraystretch}{1.0}
\setlength{\tabcolsep}{0.25in}
\begin{tabular}{c|c cc}
\cline{1-4}
$\alpha$ &  $E^{(5)}$               &$E^{(8)}$               & $E^{(10)}$ \\
\cline{1-4}
0.25    &1.50942$\times 10^{-5}$  &5.61969$\times 10^{-8}$   &1.12729$\times 10^{-9}$ \\
0.5     &1.58204$\times 10^{-5}$  &6.53078$\times 10^{-8}$   &1.32773$\times 10^{-9}$ \\
0.75    &2.91795$\times 10^{-5}$  &6.92627$\times 10^{-8}$   &1.38617$\times 10^{-9}$ \\
\hline
\end{tabular}
\label{tab2}
\end{table}

\begin{example}\label{sec3:ex2}
{\rm Consider the nonlinear singular boundary value problem
\begin{eqnarray}\label{prob_2}
\left.
  \begin{array}{ll}
(x^{\alpha}y')'=  x^{\alpha-1}e^{y}\big( -x y'-\alpha ),~~~~0<x\leq 1,\\
y(0)=\ln\frac{1}{2},~~y(1)=\ln\frac{1}{3},
\end{array}
\right\}
\end{eqnarray}
and exact solution is $y(x)=\ln(\frac{1}{2+x})$.}
\end{example}
Applying \eqref{f2} to the equation to the equation \eqref{prob_2}, with $\alpha_1=1$, $\beta_1=0$, $\gamma_1=-\ln3$ and $\eta_1=-\ln2$, we obtain
\begin{align}\label{sec3:s2}
\left.
  \begin{array}{ll}
 y_0=-\ln2,\\
 y_1=(\ln2-\ln3)x^{1-\alpha}+x^{1-\alpha}\int\limits_{0}^{1}x^{-\alpha}\int\limits_{x}^{1}x^{\alpha -1}A_{0}dxdx -\int\limits_{0}^{x}x^{-\alpha}\int\limits_{x}^{1}x^{\alpha -1}A_{0}dxdx,\\
 y_{n+1}=x^{1-\alpha}\int\limits_{0}^{1}x^{-\alpha}\int\limits_{x}^{1}x^{\alpha- 1}A_{n}dxdx-
\int\limits_{0}^{x}x^{-\alpha}\int\limits_{x}^{1}x^{\alpha -1}A_{n}dxdx,~~n\geq1.
  \end{array}
\right\}
\end{align}
The Adomian polynomials   for  $f(x,y,y')=-\big( x e^{y} y'+e^{y}\alpha )$ with $y_0=-\ln(2)$ are obtained as
\begin{align}\label{sec3:a2}
\left.
 \begin{array}{ll}
 A_0=- e^{y_0}\big(x y'_0+\alpha\big)  \\
 A_1=- e^{y_0}\big(x y'_1+y_1\alpha  \big),\\
 A_2=-e^{y_0}\big(x y'_2+y_1 y'_1)+( y_2+\nicefrac{y_1^2}{2} )\alpha \big),\\
 A_3=- e^{y_0}\big(x \left( y'_3+y'_2y_1+y'_1\left(\nicefrac{y_1^2}{2}+y_2\right)\right)+\left( y_3+y_1 y_2 +\nicefrac{y_1^3}{6} \right)\alpha  \big),\\
\ldots
\end{array}
\right\}
\end{align}
For~$\alpha=0.5$:  Using \eqref{sec3:s2} and \eqref{sec3:a2},  we obtain 
\begin{align*}
y_0&=-0.693147,\\
y_1&=0.0945349x^{0.5}-0.5x,\\
y_2&=-0.0934884 x^{0.5}-0.0315116 x^{1.5}+0.125 x^{2},\\
y_3&=-0.00413483 x^{0.5}+0.0311628 x^{1.5}-0.00111711 x^{2}+0.0157558 x^{2.5}-0.0416667 x^{3},\\
y_4&=0.00321966 x^{0.5}+0.00137828 x^{1.5}+0.00220948 x^{2}-0.0156096 x^{2.5}+0.00105504 x^{3}\\&-0.00787791 x^{3.5}+0.015625 x^{4},\\
\ldots
\end{align*}
For different values of $\alpha$, the maximum absolute error $E^{(n)}$, for $n=5,8$ and $10$, are given in  Table \ref{tab3}.  It can be noted that  when $n$  increases the maximum error decreases.
\begin{table}[htbp]\caption{Maximum error of Example~\ref{sec3:ex2}}
\centering
\vspace{-.3cm}
\renewcommand{\arraystretch}{1.0}
\setlength{\tabcolsep}{0.25in}
\begin{tabular}{c|c cc}
\cline{1-4}
$\alpha$ &  $E^{(5)}$               &$E^{(8)}$               & $E^{(10)}$ \\
\cline{1-4}
0.25    &8.25847$\times 10^{-5}$  &8.22973$\times 10^{-7}$   &3.94246$\times 10^{-8}$ \\
0.5     &8.41014$\times 10^{-5}$  &1.68118$\times 10^{-6}$   &4.29665$\times 10^{-8}$ \\
0.75    &3.17953$\times 10^{-5}$  &8.27238$\times 10^{-6}$   &6.58700$\times 10^{-8}$ \\
\hline
\end{tabular}
\label{tab3}
\end{table}

\begin{example}\label{sec3:ex3}
{\rm Consider the  singular boundary value  problem
\begin{eqnarray}\label{prob_3}
\left.
  \begin{array}{ll}
    (x^{\alpha}y')'=\beta x^{\alpha+\beta-2}\left( xy'+y(\alpha+\beta-1)\right),~~0<x\leq1, \\
  y(0)=1~~~ {\rm and}~~~ y(1)=e,
\end{array}
\right\}
\end{eqnarray}
with exact solution $y(x)=e^{x^\beta}$. }
\end{example}
Applying \eqref{f2} to the equation to the equation \eqref{prob_3}, with  $\alpha_1=1$, $\beta_1=0$, $\gamma_1=e$ and $d=\eta_1$, we obtain
\begin{align}\label{sec3:s3}
\left.
  \begin{array}{ll}
    y_0=1,\\
y_1=(e-1)x^{1-\alpha}+x^{1-\alpha}\int\limits_{0}^{1}x^{-\alpha}\int\limits_{x}^{1}x^{\alpha +\beta -2}A_{0}dx dx -
\int\limits_{0}^{x}x^{-\alpha}\int\limits_{x}^{1}x^{\alpha +\beta -2}A_{0}dxdx,\\
 y_{n+1}=x^{1-\alpha}\int\limits_{0}^{1}x^{-\alpha}\int\limits_{x}^{1}x^{\alpha +\beta -2}A_{n}dxdx-
\int\limits_{0}^{x}x^{-\alpha}\int\limits_{x}^{1}x^{\alpha +\beta -2}A_{n}dxdx,~~n\geq1,
  \end{array}
\right\}
\end{align}
where $A_n=\beta  \left( xy'_n+y_n(\alpha+\beta-1)\right).$

For~$\alpha=0.5$,~$\beta=1$: Using \eqref{sec3:s3},  we obtain
\begin{align*}
y_0&=1,\\
y_1&=0.718282 x^{0.5}+x,\\
y_2&=-0.978855 x^{0.5}+0.478855 x^{1.5}+0.5 x^{2},\\
y_3&=0.294361 x^{0.5}-0.65257 x^{1.5}+0.191542 x^{2.5}+0.166667 x^{3},\\
y_4&=-0.0316058 x^{0.5}+0.196241 x^{1.5}-0.261028 x^{2.5}+0.0547262 x^{3.5}+0.0416667 x^{4},\\
\ldots
\end{align*}
 
For $\alpha=0.5$,~$\beta=2.5$:  Using \eqref{sec3:s3},  we obtain
\begin{align*}
y_0&=1,\\
y_1&=0.718282 x^{0.5}+ x^{2.5},\\
y_2&=-1.09857 x^{0.5}+0.598568 x^{3}+0.5 x^{5},\\
y_3&=0.47673 x^{0.5}-0.915473 x^{3}+0.272076 x^{5.5}+0.166667 x^{7.5},\\
y_4&=-0.107842 x^{0.5}+0.397275 x^{3}-0.416124 x^{5.5}+0.0850239 x^{8}+0.0416667 x^{10},\\
\ldots
\end{align*}
The maximum absolute error for different values of $\alpha$ and $\beta$ in Table\ref{tab4} and  Table \ref{tab5}.
 \begin{table}[htbp]\caption{Maximum error of Example~\ref{sec3:ex3},when $\beta=1$}
\centering
\vspace{-.3cm}
\renewcommand{\arraystretch}{1.0}
\setlength{\tabcolsep}{0.25in}
\begin{tabular}{c|c cc}
\cline{1-4}
$\alpha$ &  $E^{(5)}$               &$E^{(8)}$               & $E^{(10)}$ \\
\cline{1-4}
0.25    &9.47534$\times 10^{-5}$  &9.45679$\times 10^{-7}$   &2.28700$\times 10^{-9}$ \\
0.5     &6.54946$\times 10^{-5}$  &6.83854$\times 10^{-7}$   &1.04295$\times 10^{-8}$ \\
0.75    &8.17721$\times 10^{-5}$  &1.08101$\times 10^{-6}$   &1.18947$\times 10^{-8}$ \\
\hline
\end{tabular}
\label{tab4}
\end{table}
 \begin{table}[htbp]\caption{Maximum error of Example~\ref{sec3:ex3},when $\beta=2.5$}
\centering
\vspace{-.3cm}
\renewcommand{\arraystretch}{1.0}
\setlength{\tabcolsep}{0.25in}
\begin{tabular}{c|c cc}
\cline{1-4}
$\alpha$ &  $E^{(5)}$               &$E^{(8)}$               & $E^{(10)}$ \\
\cline{1-4}
0.25    &8.48845$\times 10^{-4}$  &5.11258$\times 10^{-6}$   &3.01254$\times 10^{-8}$ \\
0.5     &9.99777$\times 10^{-4}$  &5.94699$\times 10^{-6}$   &3.44033$\times 10^{-8}$ \\
0.75    &2.97716$\times 10^{-4}$  &1.53819$\times 10^{-6}$   &1.19534$\times 10^{-8}$ \\
\hline
\end{tabular}
\label{tab5}
\end{table}

\section{Conclusion}
 The Adomian decomposition method has been used  for  solving a class of nonlinear singular boundary value problems  with derivative dependent  source function, i.e., $f(x,y,y')$.  The main advantage of the approach is that it provides a direct scheme for solving the doubly singular boundary value  problems.  The method provides a reliable technique which requires less work compared to standard Adomian decomposition method. The numerical results of the examples are presented  and only a few terms are required to obtain accurate solutions.  By comparing the results with other existing methods, it has been proved that  ADM is a more powerful method for solving doubly  singular boundary value problems.


\end{document}